\documentclass[12pt]{article}

\usepackage{latexsym}
\usepackage{amssymb}

\newcommand{\ga}{\alpha}

\newcommand{\gd}{\delta}

\newcommand{\gk}{\kappa}
\newcommand{\gl}{\lambda}

\newcommand{\rest}{\restriction}
\newcommand{\la}{\langle}
\newcommand{\ra}{\rangle}

\newcommand{\ha}{\aleph}

\newcommand{\card}[1]{{\vert #1 \vert} }
\newcommand{\ot}[1]{\hbox{o.t.($#1$)}}
\newcommand{\forces}{\Vdash}

\renewcommand{\models}{\vDash}

\newcommand{\dom}{{\rm dom}}

\newcommand{\FP}{{\mathbb P}}
\newcommand{\FQ}{{\mathbb Q}}

\newtheorem{theorem}{Theorem}

\newtheorem{conjecture}{Conjecture}
\newenvironment{proof}{\noindent{\bf
Proof:}}{\nopagebreak\mbox{}\newline\makebox[\textwidth]{\hfill$\square$}
\par\bigskip}
\newenvironment{sketch}{\noindent{\bf
 Sketch of Proof:}}{\nopagebreak\mbox{}\newline
 \makebox[\textwidth]{\hfill$\square$}\par\bigskip}

\def\tlt{\triangleleft}
\def\k{\kappa}
\def\a{\alpha}
\def\b{\beta}
\def\d{\delta}
\def\s{\sigma}
\def\t{\tau}
\def\l{\lambda}
\def\lted{{{\leq}\d}}
\def\ltk{{{<}\k}}
\def\P{{\mathbb P}}
\def\Q{{\mathbb Q}}
\def\Qdot{\dot\Q}
\def\Pforces{\forces_{\P}}

\def\card#1{\left|#1\right|}
\def\boolval#1{\mathopen{\lbrack\!\lbrack}\,#1\,\mathclose{\rbrack\!
        \rbrack}}

\def\st{\mid}
\def\set#1{\{\,{#1}\,\}}

\def\muchgt{>>}
\def\cof{\mathop{\rm cof}\nolimits}
\def\iff{\mathrel{\leftrightarrow}}
\def\intersect{\cap}
\def\minus{\setminus}
\def\Union{\bigcup}

\def\and{\mathrel{\kern1pt\&\kern1pt}}
\def\image{\mathbin{\hbox{\tt\char'42}}}
\def\elesub{\prec}
\def\iso{\cong}
\def\<#1>{\langle\,#1\,\rangle}
\def\ot{\mathop{\rm ot}\nolimits}
\title{Indestructible Weakly Compact Cardinals and the
Necessity of Supercompactness for Certain Proof Schemata
\thanks{Both authors wish to thank the thank the CUNY Research
Foundation for the Collaborative Incentive Grant, as well
as PSC-CUNY grants, which partially supported their research.}}
\author{Arthur W.~Apter\\
        Department of Mathematics\\
        Baruch College of CUNY\\
        New York, New York 10010\\
        http://math.baruch.cuny.edu/$\sim$apter\\
        awabb@cunyvm.cuny.edu\\
        \\
        Joel David Hamkins\\
        Department of Mathematics\\
        CUNY, College of Staten Island\\
        Staten Island, New York 10314\\
        http://www.math.csi.cuny.edu/$\sim$hamkins\\
        hamkins@math.csi.cuny.edu}

\date{July 4, 1999}

\begin{document}

\maketitle

\begin{abstract} We show that if the weak compactness of a cardinal
is made indestructible by means of any preparatory forcing of a
certain general type, including any forcing naively resembling the
Laver preparation, then the cardinal was originally supercompact.
We then apply this theorem to show that the hypothesis of
supercompactness is necessary for certain proof schemata.

\end{abstract}
\baselineskip=24pt

\section{Introduction and Preliminaries}\label{s1}

The well-known Laver preparation \cite{L} of a supercompact cardinal
is one of the most fundamental tools in the field of large cardinals
and forcing, and forms the backbone of many relative consistency
proofs. The first author, for example, used the Laver preparation in
\cite{A83} and \cite{A85} to produce models of ZF having long
sequences $\d$, $\d^{+}$, $\d^{++},\ldots$ of consecutive Ramsey
cardinals. While his proofs needed supercompactness in order to
apply the Laver preparation, they ultimately made use only of
indestructible weakly compact or Ramsey cardinals, an apparently
much weaker hypothesis. At the time, therefore, the
supercompactness hypotheses seemed overly strong.

The purpose of this paper is to show, however, that the
supercompactness hypotheses are far from superfluous and are
actually necessary in proofs of these theorems that employ
indestructible weakly compact cardinals produced by partial
orderings of a certain general type. Specifically, our Main Theorem,
Theorem \ref{ht0} in Section \ref{s2}, shows that after any closure
point forcing (a broad class of forcing notions defined in Section
\ref{s2}), if the weak compactness of a cardinal $\k$ becomes
indestructible by $\ltk$-directed closed forcing (originally called
``$\k$-directed closed" by Laver in \cite{L}), then $\gk$ was initially
supercompact. We then apply this in Section \ref{s3} to show that
the hypothesis of supercompactness in the proof schemata of
\cite{A83} and \cite{A85} cannot be reduced.

The terminology we use in this paper is fairly standard. For
anything left unexplained, the reader is urged to consult either \cite{J}
or \cite{K}.

\section{Indestructible Weakly Compact Cardinals}\label{s2}

Given that the Laver preparation makes any supercompact cardinal
$\k$ indestructible by $\ltk$-directed closed forcing, one naturally
wonders whether other kinds of large cardinals can be made
indestructible. Can one force, for example, to make any weakly
compact cardinal indestructible? While one might hope for a
Laver-like preparation for weakly compact cardinals, the only known
method of producing an indestructible weakly compact cardinal is to
begin with a supercompact cardinal and apply the Laver preparation
or some similar forcing. Since the strength of this hypothesis, a
supercompact cardinal, seems at first out of proportion with the
strength of the conclusion, a mere weakly compact cardinal, one
wonders whether such a strong hypothesis is really necessary. The
main theorem of this paper is that indeed for Laver-like preparations
one must begin with a supercompact cardinal, for if a weakly
compact cardinal is made indestructible by any member of a wide
class of forcing notions, including the Laver preparation and any
iteration even naively resembling it, then it was originally
supercompact. This theorem therefore improves to the weakly
compact case a previous result of the second author, which asserted
the same fact for indestructible measurable cardinals (see Corollary
5.4 of \cite{Hama}).

The class of forcing notions for which the theorem holds is quite
general. All that is required is that it admit a {\it closure point} at
some $\d$ below the cardinal in question, meaning that it factors as
$\P*\Qdot$, where $\card{\P}\leq\d$ and $\forces_\FP``\Qdot$ is
$\lted$-strategically closed''. Such a forcing notion will be referred to
as {\it closure point forcing}, particularly when the closure point
$\delta$ is smaller than whatever other cardinal $\k$ is under
discussion. This concept generalizes the concept of gap forcing in
\cite{Ham98a}, \cite{Hamb} and \cite{Ham99}, for which it was further
required that $\card{\P}<\d$. The Laver preparation admits a closure
point between any two nontrivial stages of forcing, as does Silver
forcing (the reverse Easton iteration adding a Cohen subset at every
regular
cardinal stage) and most of the other reverse Easton iterations of
closed forcing
that one commonly finds in the literature. Since closure point forcing
is so widespread, we hope that the Main Theorem will find
a wide applicability.

\newtheorem{MainTheorem}[theorem]{Main Theorem}
\begin{MainTheorem}\label{ht0}
If after forcing with a closure point below $\k$ the
weak compactness of $\k$ becomes indestructible by $\ltk$-directed
closed forcing, then $\k$ was originally supercompact.
\end{MainTheorem}

The proof makes a critical use of a generalization of the Key Lemma
of \cite{Hamb}, which asserts that forcing with a closure point at $\d$
adds no fresh $\l$-sequences of ordinals for $\cof(\l)>\d$, meaning
that if all the initial segments of a $\l$-sequence of ordinals are in the
ground model, then the sequence itself is in the ground model. Here,
we will relax the requirement that all initial segments are in the
ground model to the assumption merely that all the approximations
to the sequence with sets of size $\d$ lie in the ground model.
Specifically, if $s$ is a subset of a set in $V$, with $s$ itself not
necessarily in $V$, then the {\it $\d$-approximations} to $s$ over
$V$ are precisely the sets of the form $\s\intersect s$ for $\s\in V$ of
size $\d$ in $V$.

\newtheorem{ApproximationLemma}[definition]{Approximation
Lemma}
\begin{ApproximationLemma}\label{hl1}
Forcing with a closure point at $\d$ adds no
new subset of the ground model $V$ all of whose
$\d$-approximations are in $V$.
\end{ApproximationLemma}

\begin{proof}
By enumerating sets in $V$, it suffices to consider only sets of
ordinals. So, suppose that $V[g][H]$ has a closure point at $\d$, so
that $g*H\subseteq\P*\Qdot$ is $V$-generic, $\card{\P}\leq\d$,
$\Pforces``\Qdot$ is $\lted$-strategically closed'', and all the
$\d$-approximations to a set $s\subseteq\theta$ over $V$ are in $V$.
We want to show that $s$ itself is in $V$.

If $s$ is not in $V$, then by chopping $s$ off at an earlier ordinal if
necessary we may assume that all the initial segments of $s$ are in
$V$. It follows that $\cof(\theta)\leq\d$, for otherwise $s$ would be a
fresh sequence added by forcing with a closure point contrary to the
Key Lemma of \cite{Hamb}, mentioned just earlier. So in $V$ we
may write $\theta=\sup\set{\theta_\a\st\a<\bar\d}$, where
$\bar\d=\cof(\theta)$ and $\<\theta_\a\st\a<\bar\d>$ is a continuous,
increasing sequence of ordinals in  $V$. Since
$s\intersect\theta_\a\in V$ for all $\a<\bar\d$, it follows by the
closure of the forcing $\Qdot$ that $s\in V[g]$ and so $s=\dot s_g$
for some $\P$-name $\dot s\in V$. Let ${\cal T}_\a=
\{\,t\subseteq\theta_\a\,\st\, t\in V$ and ${\boolval{\check t=\dot
s\intersect\theta_\a}}^\P\neq 0\,\}$. These sets are uniformly
definable in $V$, and their union ${\cal T}=\Union_\a {\cal T}_\a$
forms a tree under end-extension. Since the elements of ${\cal
T}_\a$ give rise to incompatible values for $\dot s$ and therefore to
an antichain in $\P$, it must be that $\card{{\cal T}_\a}\leq\d$, and
therefore also $\card{{\cal T}}\leq\d$. Define now
$$\s=\Bigl\{\,\b<\theta\,\Bigm|\, \exists t,t'\in {\cal T}\,\,[\,
t\intersect\b=t'\intersect\b\ {\hbox{\rm and}}\ t(\b)\neq
t'(\b)\,]\,\Bigr\}.$$ In other words, $\s$ is the set of all possible
branching points for branches through the tree ${\cal T}$. Since we
already observed that ${\cal T}$ has at most $\d$ many members,
$\s$ also has size at most $\d$. Thus, $\s\intersect s$ is a
$\d$-approximation to $s$ over $V$ and hence by assumption
$\s\intersect s\in V$. Now we are nearly done. The set $s$
determines a branch through ${\cal T}$, and the information about
which direction that branch turns at every possible branching point in
${\cal T}$ is precisely contained in the set $\s\intersect s$. Using
$\s\intersect s$ as a guide in $V$, therefore, we can direct our way
through ${\cal T}$ in exactly the same way that $s$ winds through
${\cal T}$. So $s\in V$, as desired. \end{proof}

\newtheorem{FilterExtensionLemma}[definition]{Filter Extension
Lemma}
\begin{FilterExtensionLemma}\label{hl2} Suppose that $\cal
F$ is a $\k$-complete filter in $V[G]$ on a set $D$ in $V$, that
$V[G]$ admits a closure point below $\k$, and that $\cal F$
measures every subset of $D$ in $V$. Then ${\cal F} \intersect V$ is
in $V$. That is, $\cal F$ extends a measure in $V$.

\end{FilterExtensionLemma}

\begin{proof} Suppose that the forcing admits a closure point at
$\d<\k$. By the Approximation Lemma, it suffices to show that every
$\d$-approximation to ${\cal F} \intersect V$ is in $V$. So suppose
$\s\in V$ has size at most $\d$, and consider $\s\intersect ({\cal F}
\intersect V)=\s\intersect {\cal F}$. We may assume that every
member of $\s$ is a subset of $D$, since these are the only possible
members of $\s\intersect {\cal F}$. Let $\s^*$ be obtained by closing
$\s$ under complements in $D$. Since $\s^*\intersect {\cal F}$ is a
collection of at most $\d$ many sets in the filter, it follows by the
$\k$-completeness of $\cal F$ that $A=\intersect(\s^*\intersect {\cal
F})$ is in $\cal F$. In particular, $A$ is nonempty, and so we may
choose an element $a\in A$. Observe now that if $B\in\s\intersect
{\cal F}$ then $A\subseteq B$ and consequently $a\in B$.
Conversely, if $a\in B$ and $B\in\s$ then because $a\notin D\minus
B$ it follows that $A\not\subseteq D\minus B$ and so $D\minus
B\notin {\cal F}$; by the assumption that $\cal F$ measures every
set in $V$, therefore, we conclude that $B\in {\cal F}$. Thus, we have
proved for $B\in\s$ that $B\in {\cal F}\iff a\in B$. So $\s\intersect {\cal
F}$ is precisely the set of all $B\in\s$ with $a\in B$, and this is
certainly in $V$. \end{proof}

Let us now move to the proof of the Main Theorem. In fact, we will not
need to know that $\k$ is fully indestructible, but rather only that it is
indestructible by the canonical forcing to collapse cardinals above
$\k$ down to $\k$:

\begin{theorem}\label{ht1}

If after forcing with a closure point below $\k$ the weak
compactness of $\k$ becomes indestructible by the forcing to
collapse cardinals to $\k$, then $\k$ was supercompact in the
ground model.

\end{theorem}

The canonical forcing to collapse a cardinal $\theta\geq\k$ down to
$\k$ is $\ltk$-directed closed, and so by combining the preparatory
forcing with this collapsing forcing, we obtain an iteration admitting
the same closure point. Thus, the previous theorem is an immediate
consequence of the following more local version:

\begin{theorem}\label{ht2}

Suppose $\k$ is weakly compact in a forcing extension that admits a
closure point below $\k$ and that collapses $(2^{\theta^{{<}\k}})^V$ to
$\k$. Then $\k$ was $\theta$-supercompact in $V$.

\end{theorem}

\begin{proof} Suppose that $\k$ is weakly compact in $V[G]$, a
forcing extension with a closure point at $\d<\k$, and that
$(2^{\theta^{{<}\k}})^V$ is collapsed to $\k$ there. We may exhibit
the gap explicitly by writing $V[G]=V[g][H]$, where $g*H \subseteq
\P \ast \dot \FQ$, $\card{\P}\leq\d$, and $\Pforces ``\dot \FQ$ is
$\lted$-strategically closed''. Let $\eta\muchgt\theta,\k,\card{\Qdot}$,
and choose $X\elesub V_\eta[g][H]$ of size $\k$, closed under
$\ltk$-sequences, containing $\k$, $\theta$, $\P*\Qdot$, $g*H$ as
well as every element of $\theta$ and every element of
$\wp(P_\k\theta)^V$, both of which by assumption have size $\k$ in
$V[G]$. The Mostowski collapse of $X$ has the form $M[g][\tilde
H]$, where $g*\tilde H\subseteq\P*\dot\Q$ is $M$-generic for forcing
with a closure point at $\d$. Furthermore, since $X$ was closed
under $\ltk$-sequences in $V[G]$, the same holds for $M[g][\tilde
H]$. And since $M[g][\tilde H]$ has size $\k$ and $\k$ is weakly
compact, there is an embedding $j:M[g][\tilde H]\to N[g][j(\tilde H)]$
with critical point $\k$.

Consider now the set $j\image\theta$. Since $\theta$ was collapsed
to $\k$ in $V[G]$, it follows that there is a relation $\tlt$ in $X$ on
$\k$ such that $\<\k,{\tlt}>\iso\<\theta,{\in}>$. And since $\tlt$ is
fixed by the Mostowski collapse of $X$, it follows that it is also in
$M[g][\tilde H]$. By the elementarity of $j$, if $\a<\k$ is the $\b^{\rm
th}$ element with respect to the relation $\tlt$ on $\k$, then
$j(\a)=\a$ is the $j(\b)^{\rm th}$ element with respect to the relation
$j(\tlt)$ on $j(\k)$. Thus, $j\image\theta$ is precisely the set
$\set{\ot_{j(\tlt)}(\a)\st \a<\k}$ in $N[g][j(\tilde H)]$. So
$j\image\theta\in N[g][j(\tilde H)]$.

Let us now show that $j\image\theta\in N$. Since the forcing
$g*j(\tilde H)$ admits a closure point at $\d$, it suffices by the
Approximation Lemma to show that every $\d$ approximation to
$j\image\theta$ over $N$ is actually in $N$. So suppose that $\s\in
N$ has size $\d$, and consider $\s\intersect (j\image\theta)$. This
set has size $\d$, so it must have the form $j\image b$ for some
$b\subseteq\theta$. By the closure of $M[g][\tilde H]$ we know that
$b\in M[g][\tilde H]$; further, since the $\tilde H$ forcing is
$\lted$-strategically closed, it must really be that $b\in M[g]$.
Since the $g$ forcing has size $\d$, we conclude that $b\subseteq
c$ for some $c\in M$ of size $\d$. Thus, $$\s\intersect
(j\image\theta)= \s\intersect (j\image b)\subseteq\s\intersect (j\image
c)\subseteq\s\intersect (j\image\theta).$$ Consequently,
$\s\intersect (j\image\theta)=\s\intersect (j\image c)$. But $c$ has
size $\d$ and so we have $j\image c=j(c)$. Furthermore, since $c\in
M$ we know $j(c)\in N$. Thus, $\s\intersect
(j\image\theta)=\s\intersect j(c)$ is also in $N$, as we had desired.
So by the Approximation Lemma, $j\image\theta\in N$.

Let ${\cal F}=\set{X\subseteq(P_\k\theta)^V\,\st\, j\image\theta\in
j(X)}$. This is a pre-filter on $(P_\k\theta)^V$ which measures every
set in $M[g][\tilde H]$. In particular, since by design $M[g][\tilde H]$
includes every element of $\wp(P_\k\theta)^V$, we know that $\cal
F$ measures every subset of $P_\k\theta$ in $V$. And since
$M[g][\tilde H]$ is closed under $\ltk$-sequences in $V[G]$, it
follows that the filter generated by $\cal F$ is $\k$-complete in
$V[G]$. Thus, by the Filter Extension Lemma, the set $\mu={\cal
F}\intersect V$ must be in $V$.

It remains only to check that $\mu$ is a normal fine measure on
$P_\k\theta$ in $V$. Note that because $j\image\theta\in N$, we
know that $(P_\k\theta)^V\in\mu$. Certainly $\mu$ is a
$\k$-complete measure on $P_\k\theta$ in $V$, because $\cal F$ is
$\k$-complete in $V[G]$ and measures every set in $V$. It is a fine
measure because $j\image\theta$ contains the element $j(\a)$ for
every $\a<\theta$, and so for any given such $\a$, the filter $\cal F$
concentrates on the set of $\t\in P_\k\theta$ containing $\a$. To see
that $\mu$ is normal, suppose that $f:P_\k\theta\to\theta$ is
regressive in $V$. Thus, $j(f)(j\image\theta)\in j\image\theta$, and so
$j(f)(j\image\theta)=j(\a)$ for some $\a<\theta$. So $f(\t)=\a$ for
$\mu$-almost every $\t$. So $\mu$ is a normal fine measure on
$P_\k\theta$ in $V$. We conclude that $\k$ is
$\theta$-supercompact in $V$, as the theorem asserts. \end{proof}

The Main Theorem does not show that the existence of an
indestructible weakly compact cardinal is equiconsistent with the
existence of a supercompact cardinal. Rather, it shows that there is
no Laver-like preparation that makes a weakly compact
non-supercompact cardinal indestructible.
If one hopes to construct a
preparation making an arbitrary weakly compact cardinal
indestructible, therefore, one needs a completely new idea. What is
more, core model theorists \cite{Philip} report that the consistency
strength of an indestructible weakly compact cardinal is quite high,
above the upper reaches of the core model theory itself. In light of
this, and with the support suggested by our Main Theorem here, we
make the following conjecture:

\begin{conjecture}

The existence of an indestructible weakly compact cardinal is
equiconsistent
over $ZFC$ with the existence of a supercompact cardinal.
\end{conjecture}

\section{The Necessity of Supercompactness for
Certain Proof Schemata}\label{s3}

In this section, we turn our attention to showing the necessity of
supercompactness for the proof schemata given in \cite{A83} and
\cite{A85}. In order to do this, we sketch the proof of a weaker
version of Theorem 1 of \cite{A83}, illustrating the crucial use of
closure point forcing to obtain indestructibility. The Main Theorem,
therefore, implies that any modified version of this proof still
employing this crucial technique requires the use of
supercompactness. We then expand upon this to show the
necessity of the appropriate amount of supercompactness for the
proofs of Theorems 1 and 2 of \cite{A83} and Theorem 1 of \cite{A85}.

We begin by sketching the proof of
the following theorem.

\begin{theorem}\label{t1}

Let $V \models ``$ZFC + $\la \gk_i : i < \omega \ra$ is a sequence of
supercompact cardinals''. There is then a partial ordering $\FP \in V$
and a symmetric inner model $N \subseteq V^\FP$ so that $N
\models ``$ZF + $DC_{\gk_0}$ + For each $i \in (0, \omega)$,
$\gk_{i + 1}
= \gk^+_i$ + For every $i < \omega$, $\gk_i$ is a Ramsey
cardinal''.

\end{theorem}

\begin{sketch} Let $\FQ$ be the iteration of either \cite{A83} or
\cite{A98} which makes the supercompactness of each $\gk_i$
indestructible under ${<} \gk_i$-directed closed partial
orderings. Without loss of generality, as in either \cite{A83} or
\cite{A98}, we can assume that $V$ is a model constructed via
forcing over an earlier model with $\FQ$, i.e., that $V \models
``$The supercompactness of each $\gk_i$ is indestructible under
${<} \gk_i$-directed closed forcing''. We note that $\FQ$,
being an iteration of partial orderings defined in the style of \cite{L},
admits by construction a closure point below $\gk_0$.

For each $i$, let $\FP_i = {\rm Coll}(\gk_i, {<} \gk_{i + 1})$, i.e., the
usual L\'evy collapse of all cardinals in the interval $(\gk_i, \gk_{i +
1})$ to $\gk_i$. We then define $\FP = \prod_{i < \omega} \FP_i$.

Let $G$ be $V$-generic over $\FP$. $V[G]$, being a model of AC, is
not our desired model $N$. In order to define $N$, we first note that
by the Product Lemma, $G_i$, the projection of $G$ onto $\FP_i$, is
$V$-generic over $\FP_i$, and $G^i = \prod_{j < i} G_i$ is
$V$-generic over $\prod_{j < i} \FP_j$. Next, let ${\cal F} = \{f \st f :
\omega \to \sup_{i < \omega}(\k_i)$ is a function for which $f(i) \in
(\gk_i, \gk_{i + 1})\}$. For any $f \in {\cal F}$, define
$G \rest f = \{ \la p_i : i < \omega \ra\in \prod_{i < \omega} G_i \st $
For each $i < \omega$, $\dom(p_i) \subseteq \gk_i \times f(i)\}$. Note
that by the Product Lemma and the properties of the L\'evy collapse,
$G \rest f$ is $V$-generic over $\prod_{i < \omega} (\FP_i \rest f(i))$,
where $\FP_ i \rest f(i) = \{p \in \FP_i \st \dom(p) \subseteq \gk_i
\times f(i)\}$. $N$ can now intuitively be described as the least model
of ZF extending $V$ which contains, for every $f \in {\cal F}$, the set
$G \rest f$.

In order to define $N$ more formally, we let ${\cal L}_1$ be the
ramified sublanguage of the forcing language $\cal L$ with respect to
$\FP$ which contains symbols $\dot v$ for each $v \in V$, a unary
predicate symbol $\dot V$ (to be interpreted $\dot V(\dot v) \iff v \in
V$), and symbols $\dot G \rest f$ for each $f \in {\cal F}$. $N$ is
then defined as follows.

\setlength{\parindent}{1.5in}

$N_0 = \emptyset$.

$N_\gl = \bigcup_{\ga < \gl} N_\ga$ if $\ga$ is a limit ordinal.

$N_{\ga + 1} = \Bigl\{\,x \subseteq N_\ga\, \Bigm|\,
\raise6pt\vtop{\baselineskip=10pt\hbox{$x$ is definable over the
model ${\la N_\ga, \in, c \ra}_{c \in N_\ga}$} \hbox{via a term $\tau \in
{\cal L}_1$ of rank $\le \ga$}}\,\Bigr\}$.

$N = \bigcup_{\ga \in {\rm Ord}^V} N_\ga$.

\setlength{\parindent}{1.5em}

The standard arguments show $N \models {\rm ZF}$. By the proofs
of Lemmas 1.1, 1.3, 1.4, and 1.7 of \cite{A83}, $N \models ``{\rm
DC}_{\gk_0}$ + For each $i \in (0, \omega)$,
$\gk_{i + 1} = \gk^+_i$''.

To see that $N \models ``$For every $i$, $\gk_i$ is a Ramsey
cardinal'', we recapitulate the argument given in Lemma 1.2 of
\cite{A83}. Fix $i \in \omega$, and let $g : {[\gk_i]}^{< \omega} \to 2$,
$g \in N$. By Lemma 1.1 of \cite{A83}, $g \in V[\prod_{i < \omega}
G_i \rest f]$ for some $f \in {\cal F}$.

Write now $G \rest f = \prod_{j \ge i} G_j \times G^i = H \times G^i$.
By the facts that $\prod_{j \ge i} \FP_j$ is
${<} \gk_i$-directed closed,
$V \models ``\gk_i$ is an indestructible
supercompact cardinal'', and the Product Lemma, $V[H] \models
``\gk_i$ is supercompact and hence is a Ramsey cardinal''. Since in
both $V$ and $V[H]$, $|\prod_{j < i} \FP_j| < \gk_i$ (if $i = 0$, we
take this product as being the trivial partial ordering $\{\emptyset\}$),
by the L\'evy-Solovay results \cite{LS}, $V[H][G^i] = V[G \rest f]
\models ``\gk_i$ is supercompact and hence is a Ramsey cardinal''.
Thus, let $A \in V[G \rest f]$ be homogeneous for $g$. Since $V[G
\rest f] \subseteq N$, $A \in N$. This shows that $N \models ``\gk_i$
is a Ramsey cardinal''. This completes the proof sketch of Theorem
\ref{t1}.
\end{sketch}

We note that in the proof of Theorem \ref{t1}, the indestructibility of
the supercompactness of each $\gk_i$ was not used in its full force.
Rather, what was used is that each $\gk_i$ is an indestructible
Ramsey cardinal. Prima facie, this lends credence to the hope that
the hypotheses used in the proof of Theorem \ref{t1}, an $\omega$
sequence of supercompact cardinals, can be reduced in consistency
strength. Since the proof, however, uses closure point forcing to
produce the $\omega$ sequence of indestructible Ramsey cardinals,
the Main Theorem implies that for proofs still employing this
technique, the hypothesis cannot be reduced. We summarize this in
the following theorem.

\begin{theorem}\label{t2}

In any modified proof of Theorem \ref{t1} still employing closure point
forcing to produce $\omega$ many indestructible Ramsey cardinals,
the hypothesis of $\omega$ many supercompact cardinals cannot be
reduced.

\end{theorem}

Theorem 1 of \cite{A83}, in its full force, states that starting with a
model containing a proper class of supercompact cardinals, it is
possible to force and produce a model $N$ for the theory ``ZF + DC
+ Every successor cardinal is regular + The successor of every
regular cardinal is Ramsey + Every singular limit cardinal is a
Jonsson cardinal''. (In \cite{A83}, it is only mentioned that $N
\models ``$The successor of every regular cardinal is weakly
compact'', but the proofs given in Lemma 1.2 of \cite{A83} and
Theorem \ref{t1} actually show that $N \models ``$The successor of
every regular cardinal is Ramsey''.) As in the proof of Theorem
\ref{t1}, what is used to construct the relevant inner model $N$ is a
proper class of indestructible Ramsey cardinals. Since this proper
class of indestructible Ramsey cardinals is produced via a proper
class iteration of the forcing of \cite{L}, and hence admits a closure
point below the least supercompact cardinal, the Main Theorem has
the following consequence.

\begin{theorem}\label{t3}

In any modified proof of Theorem 1 in \cite{A83} still employing
closure point forcing to produce a proper class of indestructible
Ramsey cardinals, the hypothesis of a proper class of supercompact
cardinals cannot be reduced.

\end{theorem}

Theorem 2 of \cite{A83} begins with $\gk < \gd < \gl$, where $\gk$ is
a supercompact limit of supercompact cardinals, $\gd$ is
$\gl$-supercompact,
and $\gl$ is measurable, to produce a model for the
theory ``ZF + For every $i \in (0, \omega)$, $\ha_i$ is a Ramsey
cardinal + $\ha_\omega$ is a Rowbottom cardinal carrying a
Rowbottom filter + $\ha_{\omega + 1}$ is a Ramsey cardinal +
$\ha_{\omega + 2}$ is a measurable cardinal''. A key aspect of this
construction involves preliminary forcing to transform $\gk$ into a
supercompact limit of indestructible Ramsey cardinals and $\gd$ into
a cardinal whose Ramseyness is indestructible by forcing with
${<} \gd$-directed closed partial
orderings of rank less than $\gl$. Once again, since the forcing to
accomplish this is an iteration of the forcing of \cite{L}, it has
cardinality less than $\gl$ and admits a closure point below $\gk$;
and so the Main Theorem and the work of \cite{LS} yield the following
theorem.

\begin{theorem}\label{t4}

In any modified proof of Theorem 2 of \cite{A83} still producing
$\gk<\gd<\gl$, where $\gk$ is a supercompact limit of Ramsey
cardinals, $\gd$ is a Ramsey cardinal indestructible by
${<} \gd$-directed closed forcing of rank below $\gl$, and $\gl$ is
measurable, by means of forcing with a closure point below $\gk$
and of cardinality below $\gl$, the hypothesis that $\gk$ is a
supercompact limit of supercompact cardinals, that $\gd$ is
$\gl$-supercompact, and that $\gl$ is measurable cannot be
reduced.

\end{theorem}

Theorem 1 of \cite{A85} begins with an almost huge cardinal $\gk$
and disjoint sets $A, B \subseteq \gk$ for which $A \cup B = \{\ga <
\gk \st \ga$ is a successor ordinal$\}$ to construct a model $N_A$ of
height $\gk$ for the theory ``ZF + $\neg {\rm AC}_\omega$ + For
every $\ga \in A$, $\ha_\ga$ is a Ramsey cardinal + For every $\ga
\in B$, $\ha_\ga$ is a singular Rowbottom cardinal carrying a
Rowbottom filter + Every limit cardinal is a singular Jonsson cardinal
carrying a Jonsson filter''. Again, a key aspect of the proof involves a
preliminary preparation that forces $\gk$ to be a limit of Ramsey
cardinals $\gd$ that are indestructible by ${<} \d$-directed closed
forcing of rank less than $\gk$. (The almost hugeness of $\gk$ is
preserved as well.) As before, since the forcing to accomplish this is
an iteration of the forcing of \cite{L}, it admits a closure point below
$\gk$, and so the Main Theorem implies the following.

\begin{theorem}\label{t5}

In any modified proof of Theorem 1 of \cite{A85} still employing
closure point forcing to produce a cardinal $\gk$ that is a limit of
Ramsey cardinals $\gd$ indestructible by ${<}\gd$-directed closed
forcing of rank below $\gk$, the hypothesis that $\k$ is a limit of
$\ltk$-supercompact cardinals cannot be reduced.

\end{theorem}

\section{Concluding Remarks}\label{s4}

In conclusion, we would like to remark that Theorems \ref{t2} - \ref{t5}
do not show that the supercompactness hypotheses of the theorems
to which they refer are strictly required. Rather, as the title of this
paper indicates, they only make this conclusion for proofs of these
theorems having certain general features. It is conceivable that a
completely different method of proof of Theorem 1 of \cite{A83} could
be found, for example, in which the supercompactness hypothesis is
reduced. We conjecture,  and Theorems \ref{t2} - \ref{t5} suggest,
however, that the cardinal patterns given in Theorems 1 and 2 of
\cite{A83} and Theorem 1 of \cite{A85} all have consistency strength
in the realm of supercompactness or beyond.

There is additional evidence to support the conjecture that
considerable strength is required. The recent work of Schindler
\cite{S} shows, modulo certain technical assumptions necessary in
order to carry out the relevant core model arguments, that the
aforementioned cardinal patterns all have consistency strength of at
least one Woodin cardinal. Further, we have conjectured in Section
\ref{s2} that the existence of an indestructible weakly compact
cardinal is equiconsistent over ZFC with the existence of a
supercompact cardinal. In summary, since additionally all of these
cardinal patterns seem to be unattainable by forcing over a model of
AD, we expect that hypotheses minimally at the strength of
supercompactness will turn out to be necessary to produce them.


\begin{thebibliography}{99}

\bibitem{A98} A.~Apter,
``Laver Indestructibility and the Class of
Compact Cardinals'', {\it Journal of
Symbolic Logic 63}, 1998, 149--157.

\bibitem{A83} A.~Apter,
``Some Results on Consecutive Large
Cardinals'',
{\it Annals of Pure and Applied Logic 25},
1983, 1--17.

\bibitem{A85} A.~Apter,
``Some Results on Consecutive Large Cardinals
II: Applications of Radin Forcing'',
{\it Israel Journal of Mathematics 52},
1985, 273--292.

\bibitem{Ham98a} J.~D.~Hamkins,
``Destruction or Preservation As You
Like It'',
{\it Annals of Pure and Applied Logic 91},
1998, 191--229.

\bibitem{Hamb} J.~D.~Hamkins, ``Gap Forcing'',
submitted for publication to the
{\it Journal of Mathematical Logic}.

\bibitem{Ham99} J.~D.~Hamkins, ``Gap Forcing:
Generalizing the L\'evy-Solovay Theorem'',
{\it Bulletin of Symbolic Logic 5}, 1999, 264--272.

\bibitem{Hama} J.~D.~Hamkins,
``The Lottery Preparation'', to
appear in the {\it Annals of Pure and Applied Logic}.

\bibitem{Ham98b} J.~D.~Hamkins, ``Small Forcing Makes
Every Cardinal Superdestructible'',
{\it Journal of Symbolic Logic 63}, 1998, 51--58.

\bibitem{Ham98c} J.~D.~Hamkins, S.~Shelah,
``Superdestructibility: A Dual to Laver's
Indestructibility'',
{\it Journal of Symbolic Logic 63}, 1998, 549--554.

\bibitem{J} T.~Jech, {\it Set Theory},
Academic Press, New York and San
Francisco, 1978.

\bibitem{K} A.~Kanamori, {\it The
Higher Infinite}, Springer-Verlag,
Berlin and New York, 1994.

\bibitem{L} R.~Laver,
``Making the Supercompactness of $\gk$
Indestructible under $\gk$-Directed
Closed Forcing'',
{\it Israel Journal of Mathematics 29},
1978, 385--388.

\bibitem{LS} A.~L\'evy, R.~Solovay,
``Measurable Cardinals and the Continuum Hypothesis'',
{\it Israel Journal of Mathematics 5}, 1967, 234--248.

\bibitem{S} R.~Schindler,
``Successive Weakly Compact or Singular Cardinals'',
{\it Journal of Symbolic Logic 64}, 1999, 139--146.

\bibitem{Philip} P.~Welch, 1998, personal communication with the
second author.

\end{thebibliography}
\end{document}